\documentclass[draft,11pt,english]{article}

\usepackage{amsfonts,amsmath,amsthm,amscd,amssymb,latexsym,cite,verbatim,texdraw,floatflt,
caption2,pb-diagram}
\usepackage[mathscr]{eucal}

\theoremstyle{plain}
\newtheorem{theorem}{Theorem}%[section]

\newtheorem{corollary}{Corollary}[section]

\theoremstyle{definition}

\newcommand{\keywords}{\textbf{Key words. }\medskip}
\newcommand{\subjclass}{\textbf{MSC 2010. }\medskip}
\renewcommand{\abstract}{\textbf{Abstract. }\medskip}

\begin{document}

\sloppy

\title{Asymptotic estimates for the widths \\ of classes of functions of high smoothness \thanks{This work was partially supported by the Grant H2020-MSCA-RISE-2019, project number 873071 (SOMPATY: Spectral Optimization: From Mathe\-ma\-tics to Physics and Advanced Technology).%Authors want to thanks \dots
}}

\author{A.S.~Serdyuk and I.V.~Sokolenko}

%\shorttitle{Asymptotic estimates for best uniform approximations}

%\shortauthor{A.S.~Serdyuk, I.V.~Sokolenko}

\date{\ }

\maketitle

\begin{abstract}
We find two-sided estimates for Kolmogorov, Bernstein, linear and projection widths of the classes of  convolutions of $2\pi$-periodic functions $\varphi$, such that $\|\varphi\|_2\le1$, with fixed generated kernels $\Psi_{\bar{\beta}}$, which have Fourier series of the form $\sum\limits_{k=1}^\infty \psi(k)\cos(kt-\beta_k\pi/2), $ where  $\psi(k)\ge0,$ $\sum\psi^2(k)<\infty, \beta_k\in\mathbb{R},$ in the space $C$.
It is shown that for  rapidly decrising sequences $\psi(k)$ (in particular, if $\lim\limits_{k\rightarrow\infty}{\psi(k+1)}/{\psi(k)}=0$) obtained estimates are asymptotic equalities. We establish that asymptotic equalities for  widths of this classes are realized by trigonometric Fourier sums.
\end{abstract}



\section*{Introduction}
Let   $L_{p}$,
$1\le  p<\infty$, be the space of $2\pi$--periodic functions $f$ sum\-mable to the power $p$ on  $[-\pi,\pi)$, in which the norm is given by the formula
$$
\|f\|_{L_p}=\|f\|_{p}=\bigg(\int\limits_{-\pi}^{\pi}|f(t)|^pdt\bigg)^{1/p},
$$
$L_{\infty}$ be the space of measurable and essentially bounded   $2\pi$--periodic functions  $f$ with the norm
$$
\|f\|_{L_\infty}=\|f\|_{\infty}=\mathop {\rm ess \sup}\limits_{t} |f(t)|,
$$
$C$ be the space of continuous $2\pi$--periodic functions  $f$, in which the norm is defined by the equality
$$
{\|f\|_{C}=\max\limits_{t}|f(t)|}.
$$

Denote by $C^\psi_{\bar\beta,p},\ 1\le p\le\infty,$ the set of all $2\pi$-periodic functions $f$,
representable as convolution
\begin{equation}\label{25_7'}
	f(x)=\frac{a_0}{2}+\frac{1}{\pi}\int\limits_{-\pi}^{\pi}
	\varphi(x-t) \Psi_{\bar\beta}(t)dt, \ \ \ a_0\in\mathbb R, \ \ \ \varphi\in B_p^0,
\end{equation}
$$
B_p^0=\{g\in L_p:\ \|g\|_p\le1,\ g\perp1\},
$$
with a fixed generated kernel $\Psi_{\bar\beta}\in L_{p'},\  1/p+ 1/{p'}=1,\ $ the Fourier series of which has the form
\begin{equation}\label{1*}
	S[\Psi_{\bar{\beta}}](t)=\sum\limits_{k=1}^\infty \psi(k)\cos\left(kt-\frac{\beta_k\pi}2\right),\quad \beta_k\in\mathbb{R},\quad \psi(k)\ge0.
\end{equation}
A function $f$ in the representation (\ref{25_7'}) is called $(\psi,\bar{\beta})$-integral of the function $\varphi$ and  is denoted by ${\mathcal J}^{ \psi}_{\bar{\beta}}\varphi$  $(f={\mathcal J}^{ \psi}_{\bar{\beta}}\varphi)$. If $\psi(k)\neq0,\ k\in\mathbb{N},$ then the function $\varphi$ in the representation  (\ref{25_7'}) is called $(\psi,\bar{\beta})$-derivative of the function $f$ and is denoted by $f^{ \psi}_{\bar{\beta}}$ $(\varphi=f^{ \psi}_{\bar{\beta}})$. The concepts of  $(\psi,\bar{\beta})$-integral and $(\psi,\bar{\beta})$-derivative  was introduced by Stepanets (see, e.g., \cite{Stepanets1987,Stepanets2005}).
Since \mbox{$\varphi\in L_p$} and $\ \Psi_{\bar\beta}\in L_{p'},$ then  the function $f$ of the form \eqref{25_7'} is a continuous function, i.e. $C^\psi_{\bar\beta,p}\subset C$ (see \cite[Proposition 3.9.2.]{Stepanets2005}).

In the case $\beta_k\equiv\beta,\ \beta\in \mathbb R,$ the classes $C^\psi_{\bar\beta,p}$ are denoted by $C^\psi_{\beta,p}$.

For $\psi(k)=k^{-r}, r>0, $ the classes  $C^\psi_{\bar\beta,p}$ and $C^\psi_{\beta,p}$ are denoted by $W^r_{\bar\beta,p}$ and $W^r_{\beta,p}$, respectively. The classes $W^r_{\beta,p}$ are  the well-known  Weyl-Nagy classes (see, e.g., \cite{Sz.-Nagy1938, Stechkin1956_2, Stepanets1987, Stepanets2005}). In other words, $W^r_{\beta,p}, 1\le p\le \infty,$ are the classes of  $2\pi$-periodic functions $f$, representable
as convolutions of the form
\begin{equation}\label{1}
	f(x)=\frac{a_0}{2}+\frac{1}{\pi}\int\limits_{-\pi}^{\pi}
	\varphi(x-t) B_{r,\beta}(t)dt, \ \ \ a_0\in\mathbb R,
\end{equation}
the Weyl-Nagy kernels $B_{r,\beta}$ of the form
\begin{equation}\label{2*}
	B_{r,\beta}(t)=\sum\limits_{k=1}^\infty k^{-r}\cos\left(kt-\frac{\beta\pi}2\right),\quad r>0,\quad \beta\in\mathbb{R},
\end{equation}
with functions $\varphi \in B_p^0$. The function $\varphi$ in the formula \eqref{1} is called  the Weyl-Nagy derivative of the function $f$ and is denoted by $f_\beta^r$.

If  $r\in\mathbb N$ and $ \beta=r,\ $ then the functions $B_{r,\beta}$ of the form (\ref{2*})  are  the well-known Bernoulli kernels and  the corresponding classes $W^r_{\beta,p}$ coincide with the well-known classes $W^r_{p}$ which consist of $2\pi$-periodic functions $f$  with absolutely continuous derivatives $f^{(k)}$ up to $(r-1)$-th order inclusive and such that  $\|f^{(r)}\|_p\le1$. In addition,   for almost everywhere $x\in\mathbb{R} \ \ f^{(r)}(x)=f_r^r(x)=\varphi(x),\ $  where $\varphi$ is the function from (\ref{1}).

For $\psi(k)=e^{-\alpha k^{r}},\ \alpha>0, \ r>0, $   the classes  $C^\psi_{\bar\beta,p}$ and $C^\psi_{\beta,p}$ are denoted by $C^{\alpha,r}_{\bar\beta,p}$ and $C^{\alpha,r}_{\beta,p}$, respectively.  The sets $C^{\alpha,r}_{\beta,p}$ are  well-known classes of the generalized Poisson integrals (see, e.g., \cite{Stepanets1987, Stepanets2005}), i.e. classes of  convolutions
\begin{equation}\label{'1}
	f(x)=\frac{a_0}{2}+\frac{1}{\pi}\int\limits_{-\pi}^{\pi}
	\varphi(x-t) P_{\alpha,r,\beta}(t)dt, \quad a_0\in\mathbb R,\quad \varphi \in B_p^0,
\end{equation}
with the generalized Poisson kernels
\begin{equation}\label{'2*}
	P_{\alpha,r,\beta}(t)=\sum\limits_{k=1}^\infty e^{-\alpha k^{r}} \cos\left(kt-\frac{\beta\pi}2\right),\quad \alpha>0, \quad r>0,\quad \beta\in\mathbb{R}.
\end{equation}

Let $\mathfrak N$ be a some functional class from the space $C$ $(\mathfrak N\subset C)$. The quantity
\begin{equation}\label{15_3'}
	E_n(\mathfrak N)_C=\sup\limits_{f\in \mathfrak N} E_n(f)_C= \sup\limits_{f\in \mathfrak N}\inf\limits_{T_{n-1}\in{\mathcal T}_{2n-1}} \|f-T_{n-1}\|_C
\end{equation}
is called the best uniform approximation of the class $\mathfrak{N}$ by elements of the
subspace ${\mathcal T}_{2n-1}$ of trigonometric polynomials $T_{n-1}$ of the order $n-1$:
$$
T_{n-1}(x)=\frac{\alpha_0}{2}+\sum\limits_{k=1}^{n-1}
(\alpha_k\cos kx+\beta_k\sin kx),\ \ \  \alpha_k,\beta_k\in \mathbb{R}.
$$

The order estimates for the best approximations $E_n(K)_C$ of classes \mbox{$K=C^\psi_{\bar\beta,p},$} $ 1\le p\le\infty,$ (and, hence, classes $W^r_{\beta,p}$, $C^{\alpha,r}_{\beta,p}$ and $C^\psi_{\beta,p}$) depending on rate of decreasing to zero of sequences $\psi(k)$  were obtained, in particular, in the works of Temlyakov \cite{Temlyakov1993}, Hrabova and Serdyuk \cite{Hrabova_Serdyuk2013}, Serdyuk and Stepanyuk \cite{Serdyuk_Stepanyuk2014_UMZh9,Serdyuk_Stepanyuk2014_UMZh12}.% in the works \cite{Temlyakov1993,Hrabova_Serdyuk2013,Serdyuk_Stepanyuk2014_UMZh9,Serdyuk_Stepanyuk2014_UMZh12}.

If the sequences $\psi(k)$ decrease to zero faster than any geometric progression, then asymptotic equations of the best uniform approximations are even known (see, for example, the authors work \cite{Serdyuk_Sokolenko2020} and the bibliography available there).

In \cite{Serdyuk_Sokolenko2020} it was shown that for such classes $C^\psi_{\bar\beta,p}$ the following asymptotic equations take places
\begin{equation}\label{**}
	E_n(C^\psi_{\bar\beta,p})_C\sim {\mathcal E}_n(C^\psi_{\bar\beta,p})_C\sim\frac{\|\cos t\|_{p'}}{\pi}\psi(n), \quad1\le p\le\infty,	
\end{equation}
where
\begin{equation*}\label{key}
	{\mathcal E}_n(C^\psi_{\bar\beta,p})=\sup_{f\in C^\psi_{\bar\beta,p}}\|f-S_{n-1}(f)\|_C,
\end{equation*}
$S_{n{-}1}(f)$ is the partial Fourier sum of order $n{-}1$ of the function $f$, \mbox{$\displaystyle\frac1p+\frac1{p'}=1,$} and %"$\sim$"\ is asymptotically equivalent to  i.e.
$A(n){\sim} B(n)$ as $n\rightarrow\infty$ means that $\lim\limits_{n\rightarrow\infty} {A(n)}/{B(n)}{=}1.$

For $p=\infty$ in the case of $K=W^r_{\bar\beta,\infty}, r>0, $ and in the cases of $K=C^{\alpha,r}_{\bar\beta,\infty}, r\ge1,$ and  $K=C^\psi_{\bar\beta,\infty}$ $(K=C^\psi_{\beta,\infty})$  for certain restrictions on sequences $\psi$ and $\bar{\beta}$ the exact values of the best uniform approximations are known thanks to the works  of Favard
\cite{Favard_1936,Favard_1937},   Akhiezer and Krein \cite{Akhiezer_Krein1937}, Krein \cite{Krein_1938}, Nagy \cite{Sz.-Nagy1938}, Stechkin \cite{Stechkin1956_2}, Dzyadyk   \cite{Dzyadyk _1959,Dzyadyk _1974}, Sun \cite{Sun_1961}, Bushanskij \cite{Bushanskij}, Pinkus \cite{Pinkus1985}, Serdyuk \cite{Serdyuk1995,Serdyuk1998,Serdyuk1999,Serdyuk_2002_zb,Serdyuk_2005_7}    etc.%\cite{Favard_1936,Favard_1937,Sz.-Nagy1938,Akhiezer_Krein1937,Dzyadyk _1959,Dzyadyk _1974,Krein_1938,Stechkin1956_2,Sun_1961,Bushanskij,Serdyuk1995,Serdyuk1999,Serdyuk_2002_zb,Pinkus1985}.

For $p=2$ and for arbitrary $\bar{\beta}=\beta_k\in\mathbb{R}, \sum\limits_{k=1}^\infty\psi^2(k)<\infty$ the exact values for the quantity ${\mathcal E}_n(C^\psi_{\bar\beta,2})_C$ are also known  (see  \cite{Serdyuk_Sokolenko2011Zb}).

In this paper, we establish two-sided estimates of Kolmogorov, Bernstein, linear and projection widths of the classes $C^\psi_{\bar\beta,2}$ in the space $C$, which become into asymptotic equations under certain restrictions on the sequence $\psi(k)$ (in particular, if $\lim\limits_{k\rightarrow\infty}{\psi(k+1)}/{\psi(k)}=0$).

Let  $K$ be a convex centrally symmetric subset of $C$ and
let ${ B }$ be a unit ball of the space  $C$. Let also  $F_N$ be an arbitrary $N$-dimensional subspace of space $C$, $N\in {\mathbb N}$, and
$\mathscr{L}(C, F_N)$ be a set of linear operators from $C$ to $F_N$.
By    $\mathscr {P}(C, F_N)$ denote the subset of projection operators of the set ${\mathscr{L}}(C, F_N)$,  that is, the set of the operators  $A$ of linear projection onto the set $F_N$ such that $Af = f$ when $f\in F_N$.  The quantities
\begin{equation}\label{A}
	b_N(K, C)=\sup\limits _{F_{N+1}}\sup\{\varepsilon>0: \varepsilon { B }\cap F_{N+1}
	\subset K\},
\end{equation}
\begin{equation}\label{B}
	d_N(K, C)=\inf\limits _{F_N}\sup \limits _{f\in K}
	\inf \limits _{u\in F_N}\|f - u \|_{_{\scriptstyle  C}} ,
\end{equation}
\begin{equation}\label{C}
	\lambda _N(K,C)=
	\inf \limits _{F_N}\inf \limits_{A\in {\mathscr {L}}(C, F_N)}\sup \limits _{f\in K}
	\|f - Af\|_{_{\scriptstyle  C}} ,
\end{equation}
\begin{equation}\label{D}
	\pi _N(K, C)=\inf \limits_{F_N}\inf \limits _{A\in {\mathscr {P}}(C,F_N)}
	\sup \limits _{f\in K}\|f - Af\|_{_{\scriptstyle  C}} ,
\end{equation}
are called Bernstein, Kolmogorov, linear, and projection $N$-widths of the set $K$ in the space $C$, respectively.

The results containing order estimates of the widths \eqref{A}-\eqref{D} in the case of $K=C_{\bar\beta,p}^\psi$ (and, in particular, $W^r_{\beta,p}$ and $C_{\beta,p}^\psi$) can be found, for example, in the works  of Tikhomirov \cite{Tikhomirov1976}, Pinkus \cite{Pinkus1985}, Kornejchuk \cite{Kornejchuk1987}, Kushpel' \cite{Kushpel'1989_UMZh4}, Romanyuk \cite{Romanyuk2012}, Temlyakov \cite{Temlyakov1990,Temlyakov1993}  etc. 

\section{Main results}
The main result of this paper is the following statement.

\begin{theorem}\label{theorem1}
	Let $\bar{\beta}=\{\beta_k\}_{k=1}^\infty,$ $\beta_k\in\mathbb R,$ and  $\psi(k)>0$  satisfies the condition
	\begin{equation}\label{1t1}
		\sum\limits_{k=1}^\infty\psi^2(k)<\infty.
	\end{equation}
	Then for all $n\in\mathbb{N}$ the following inequalities hold
	\begin{equation}
		\frac1{\sqrt{\pi}}\left(
		\frac1{\psi^2(n)}+2\sum\limits_{k=1}^{n-1}\frac1{\psi^2(k)}\right)^{-\frac12}\le P_{2n}(C_{\bar\beta,2}^\psi,C)
		\le P_{2n-1}(C_{\bar\beta,2}^\psi,C)%\le
		\le \frac1{\sqrt{\pi}}\left(\sum\limits_{k=n}^\infty\psi^2(k)\right)^{\frac12},
		\label{1t2}
	\end{equation}
	where $P_N$ is any of the widths $b_N, d_N, \lambda_N$ or $\pi_N$.
	\\
	\indent If, in adition, $\psi(k)$ satisfies the condition
	\begin{equation}\label{1t3}
		\lim\limits_{n\rightarrow\infty}\max\left\{
		\psi(n)\left(\sum\limits_{k=1}^{n-1}\frac1{\psi^2(k)}\right)^{\frac12},
		\frac1{\psi(n)}{\left(\sum\limits_{k=n+1}^\infty\psi^2(k)\right)^{\frac12}}	\right\}=0,
	\end{equation}
	then the following asymptotic equalities hold
	\begin{equation}\label{1t4}
		\begin{split}
			\left.\begin{array}{l}
				P_{2n}(C_{\bar\beta,2}^\psi,C)\\
				P_{2n-1}(C_{\bar\beta,2}^\psi,C)
			\end{array}\right\}= & \psi(n)\Bigg(\frac1{\sqrt{\pi}}
			\\
			+& 		\left.{\mathcal O}(1) \max\left\{	\psi(n)\left(\sum\limits_{k=1}^{n-1}\frac1{\psi^2(k)}\right)^{\frac12},
			\frac1{\psi(n)}{\left(\sum\limits_{k=n+1}^\infty\psi^2(k)\right)^{\frac12}}\right\} \right),
		\end{split}
	\end{equation}
	where ${\mathcal O}(1)$ are the quantities uniformly bounded in all parameters.
	\\
	\indent The equalities \eqref{1t4} are realized by trigonometric Fourier sums $S_{n-1}(f)$.
\end{theorem}

\begin{proof} In the work \cite{Serdyuk_Sokolenko2011Zb} it was proved that if the condition \eqref{1t1} is satisfied then the following equality holds
	\begin{equation}\label{1td1}
		{\mathcal E}(C_{\bar\beta,2}^\psi;S_{n-1})_C=\frac1{\sqrt{\pi}}	\left(\sum\limits_{k=n}^\infty\psi^2(k)\right)^{\frac12},\quad\beta_k\in\mathbb{R},\quad n\in\mathbb{N}.
	\end{equation}

	Since the operator that assigns to each function $f\in C$ its partial Fourier sum is a linear projector, then by virtue \eqref{1td1}
	\begin{equation}\label{1td2}
		\pi_{2n-1}(C_{\bar\beta,2}^\psi,C)\le {\mathcal E}(C_{\bar\beta,2}^\psi;S_{n-1})_C=\frac1{\sqrt{\pi}}	\left(\sum\limits_{k=n}^\infty\psi^2(k)\right)^{\frac12},\quad\beta_k\in\mathbb{R},\quad n\in\mathbb{N}.
	\end{equation}
	
	For all  $n\in\mathbb{N}$ and $\mathfrak{N}\subset C$
	\begin{equation*}\label{key}
		P_{2n}(\mathfrak{N},C)\le
		P_{2n-1}(\mathfrak{N},C),
	\end{equation*}
	where $P_N$ is any of the widths $b_N, d_N, \lambda_N,$ and $\pi_N$, and, in addition, for all $N\in\mathbb{N}$
	\begin{equation}\label{1td3}
		b_N(\mathfrak{N}, C)\le d_N(\mathfrak{N}, C)\le\lambda_N(\mathfrak{N}, C)\le\pi_N(\mathfrak{N}, C).
	\end{equation}
	Therefore on the basis of \eqref{1td2} we obtaine an estimate from above for the widths $P_N$ in the formula  \eqref{1t2}.
	
	To obtain a required estimate from below in \eqref{1t2}  it suffices to establish that
	\begin{equation}\label{1td4}
		b_{2n}(C_{\bar\beta,2}^\psi,C)\ge \frac1{\sqrt{\pi}}\left(
		\frac1{\psi^2(n)}+2\sum\limits_{k=1}^{n-1}\frac1{\psi^2(k)}\right)^{-\frac12}.
	\end{equation}
	
	In $(2n+1)$-dimensional space ${\mathcal T}_{2n+1}$ of trigonometric polynomials $ T_n $ of order $ n $ let us consider a ball of the form
	\begin{equation}\label{1td5}
		B_{2n+1}=\left\{T_n\in{\mathcal T}_{2n+1}:\|T_n\|_C\le \frac1{\sqrt{\pi}}\left(
		\frac1{\psi^2(n)}+2\sum\limits_{k=1}^{n-1}\frac1{\psi^2(k)}\right)^{-\frac12}\right\}
	\end{equation}
	and prove the following embedding
	\begin{equation}\label{1td6}
		B_{2n+1}\subset C_{\bar\beta,2}^\psi.
	\end{equation}
	
	For any trinometric polynomial
	\begin{equation}\label{1td7}
		T_n(x)=\frac{a_0}2+\sum\limits_{k=1}^n (a_k\cos kx+b_k\sin kx)
	\end{equation}
	from the ball $B_{2n+1}$ its $(\psi,\bar{\beta})$-derivative has a form
	$$
	(T_n)_{\bar{\beta}}^\psi(x)=\sum\limits_{k=1}^n \left(\frac{a_k}{\psi(k)}\cos \left(kx+\frac{\beta_k\pi}{2}\right)+\frac{b_k}{\psi(k)}\sin \left(kx+\frac{\beta_k\pi}{2}\right)\right)
	$$
	$$
	=\sum\limits_{k=1}^n \left(\frac{a_k\cos \frac{\beta_k\pi}{2}}{\psi(k)}\cos kx-
	\frac{a_k\sin \frac{\beta_k\pi}{2}}{\psi(k)}\sin kx+
	\frac{b_k\cos\frac{\beta_k\pi}{2}}{\psi(k)}\sin kx+
	\frac{b_k\sin\frac{\beta_k\pi}{2}}{\psi(k)}\cos kx\right)
	$$
	\begin{equation}\label{1td8}
		=\sum\limits_{k=1}^n \frac1{\psi(k)}\left(\left({a_k\cos \frac{\beta_k\pi}{2}}
		+{b_k\sin\frac{\beta_k\pi}{2}}\right)\cos kx+\left(-
		{a_k\sin \frac{\beta_k\pi}{2}}+{b_k\cos\frac{\beta_k\pi}{2}}\right)\sin kx\right).
	\end{equation}
	
	By virtue of Parseval equality, from \eqref{1td8} we get
	\begin{equation}\label{1td9}
		\left\|(T_n)_{\bar{\beta}}^\psi\right\|_2
		{=}\sqrt{\pi}	\left(\sum\limits_{k=1}^n \frac1{\psi^2(k)}\left(\left({a_k\cos \frac{\beta_k\pi}{2}}
		{+}{b_k\sin\frac{\beta_k\pi}{2}}\right)^2{+}\left({-}
		{a_k\sin \frac{\beta_k\pi}{2}}{+}{b_k\cos\frac{\beta_k\pi}{2}}\right)^2\right)\right)^{\frac12}
		$$
		$$
		=\sqrt{\pi}	\left(\sum\limits_{k=1}^n \frac1{\psi^2(k)}\left(a_k^2+b_k^2\right)\right)^{\frac12}.
	\end{equation}
	
	By Parseval equality for the polynomial $ T_n $ of the form \eqref{1td7} we obtain
	\begin{equation}\label{1td10}
		\frac{a_0^2}2+\sum\limits_{k=1}^n (a_k^2+b_k^2)=\frac1\pi\int\limits_{-\pi}^{\pi}T_n^2(x) dx.
	\end{equation}
	Therefore we have a chain of inequalities
	\begin{equation}\label{1td11}
		a_k^2+b_k^2\le \sum\limits_{k=1}^n (a_k^2+b_k^2)\le \frac1\pi\int\limits_{-\pi}^{\pi}\left\|T_n\right\|_C^2 dx= 2\left\|T_n\right\|_C^2,
	\end{equation}
	and, consequently, we obtaine an estimate for $\sqrt{a_k^2+b_k^2}$ of the following form
	\begin{equation}\label{1td12}
		\sqrt{a_k^2+b_k^2}\le\sqrt2\left\|T_n\right\|_C,\quad k=\overline{1,n}.
	\end{equation}
	
	In the case of $ k = n $ this estimate  can be improved. To do this, let us consider a trigonometric polynomial
	$$
	\tau_n(x):=\frac{T_n(x)}{\sqrt{a_n^2+b_n^2}}=\frac1{\sqrt{a_n^2+b_n^2}}\left(\frac{a_0}2+\sum\limits_{k=1}^n (a_k\cos kx+b_k\sin kx)\right)
	$$
	$$
	=\frac1{\sqrt{a_n^2+b_n^2}}\left(\frac{a_0}2+\sum\limits_{k=1}^n \sqrt{a_k^2+b_k^2} (\frac{a_k}{\sqrt{a_k^2+b_k^2}}\cos kx+\frac{b_k}{\sqrt{a_k^2+b_k^2}}\sin kx)\right)
	$$
	\begin{equation}\label{1td13}
		=\frac{\rho_0}{2}+\sum\limits_{k=1}^n\rho_k\cos(kx+\theta_k),
	\end{equation}
	where
	$$
	\rho_0=\frac{a_0}{\sqrt{a_n^2+b_n^2}},\quad\rho_k=\frac{\sqrt{a_k^2+b_k^2}}{\sqrt{a_n^2+b_n^2}},\quad k=\overline{1,n-1},\quad\rho_n=1,
	$$
	and $\theta_k$ are such that
	$$
	\left\{\begin{array}{rr}
		\cos\theta_k=\displaystyle\frac{a_k}{\sqrt{a_k^2+b_k^2}}, \ k=\overline{1,n},
		\\
		\sin\theta_k=\displaystyle\frac{-b_k}{\sqrt{a_k^2+b_k^2}}, \ k=\overline{1,n}.
	\end{array}\right.
	$$
	As it follows from \cite[statement 2.9.1]{Kornejchuk1992} for all $p\in[1,\infty]$ the following inequality holds
	$$
	\|\tau_n(\cdot)\|_p
	\ge\|\cos n(\cdot)\|_p,
	$$
	and, consequently, for $p=\infty$
	\begin{equation}\label{1td14}
		\|\tau_n\|_C
		\ge1.
	\end{equation}
	
	From \eqref{1td13} and \eqref{1td14} we get
	\begin{equation}\label{1td15}
		\sqrt{a_n^2+b_n^2}=\frac{\left\|T_n\right\|_C}{\left\|\tau_n\right\|_C}\le\left\|T_n\right\|_C.
	\end{equation}
	
	Using the equations \eqref{1td9} and the estimates \eqref{1td12} and \eqref{1td15}  we have
	\begin{equation}\label{1td16}
		\left\|(T_n)_{\bar{\beta}}^\psi\right\|_2
		=\sqrt{\pi}	\left(\sum\limits_{k=1}^{n-1} \frac{a_k^2+b_k^2}{\psi^2(k)}+\frac{a_n^2+b_n^2}{\psi^2(n)}\right)^{\frac12}
		$$
		$$
		\le \sqrt{\pi}
		\left(2\sum\limits_{k=1}^{n-1} \frac{\left\|T_n\right\|^2_C}{\psi^2(k)}+\frac{\left\|T_n\right\|^2_C}{\psi^2(n)}\right)^{\frac12}=\sqrt{\pi}
		\left(\frac{1}{\psi^2(n)}+2\sum\limits_{k=1}^{n-1} \frac{1}{\psi^2(k)}\right)^{\frac12}\left\|T_n\right\|_C.
	\end{equation}
	
	Since the polynomials $T_n$  belongs to the ball $B_{2n+1}$ of the form \eqref{1td5}, from \eqref{1td16} it follows that
	\begin{equation}\label{1td17'}
		\left\|(T_n)_{\bar{\beta}}^\psi\right\|_2\le1.
	\end{equation}
	
	The embedding \eqref{1td6} is proved.
	
	The inequality \eqref{1td4} follows from the definition of the Berstein width $b_{2n}(C^\psi_{\bar\beta,2},C)$ and the embedding \eqref{1td6}. The relations \eqref{1td2}-\eqref{1td4} prove the inequalities \eqref{1t2}.

	To prove the asymptotic equations \eqref{1t4} under satisfying the condition \eqref{1t3} first of all we note  that
	\begin{equation}\label{1td17}
		\left(\sum\limits_{k=n}^\infty \psi^2(k)\right)^{\frac12}\le \psi(n)+\left(\sum\limits_{k=n+1}^\infty \psi^2(k)\right)^{\frac12}
	\end{equation}
	and
	\begin{equation}\label{1td18}
		\left(\frac1{\psi^2(n)}	+2\sum\limits_{k=1}^{n-1} \frac{1}{\psi^2(k)}\right)^{\frac12}\le \frac1{\psi(n)}+\left(2\sum\limits_{k=1}^{n-1} \frac{1}{\psi^2(k)}\right)^{\frac12}.
	\end{equation}
	
	From \eqref{1td18} we get
	\begin{equation}\label{1td19}
		\left(\frac1{\psi^2(n)}	+2\sum\limits_{k=1}^{n-1} \frac{1}{\psi^2(k)}\right)^{-\frac12}\ge \frac1{\frac1{\psi(n)}+\left(2\sum\limits_{k=1}^{n-1} \frac{1}{\psi^2(k)}\right)^{\frac12}}
		$$
		$$
		=\psi(n)-\left(\frac1{\frac1{\psi(n)}}-\frac1{\frac1{\psi(n)}+\left(2\sum\limits_{k=1}^{n-1} \frac{1}{\psi^2(k)}\right)^{\frac12}}\right)
		=
		\psi(n)-\frac{\left(2\sum\limits_{k=1}^{n-1} \frac{1}{\psi^2(k)}\right)^{\frac12}}{\frac1{\psi(n)}\left(\frac1{\psi(n)}+\left(2\sum\limits_{k=1}^{n-1} \frac{1}{\psi^2(k)}\right)^{\frac12}\right)}
		$$
		$$
		=\psi(n)\left(1-\frac{\left(2\sum\limits_{k=1}^{n-1} \frac{1}{\psi^2(k)}\right)^{\frac12}}{\frac1{\psi(n)}+\left(2\sum\limits_{k=1}^{n-1} \frac{1}{\psi^2(k)}\right)^{\frac12}}\right)
		\ge
		\psi(n)\left(1-\psi(n)\left(2\sum\limits_{k=1}^{n-1} \frac{1}{\psi^2(k)}\right)^{\frac12}\right).
	\end{equation}
	
	So, as it follows from \eqref{1t2} and \eqref{1td17}, on the one hand,
	\begin{equation}\label{1td20}
		P_{2n}(C^\psi_{\bar\beta,2},C)\le P_{2n-1}(C^\psi_{\bar\beta,2},C)
		\le
		\psi(n)\left(\frac1{\sqrt{\pi}}+\frac1{\sqrt{\pi}\psi(n)}\left(\sum\limits_{k=n+1}^\infty \psi^2(k)\right)^{\frac12}\right),
	\end{equation}
	and, on the other hand, by virtue of \eqref{1t2}
	\begin{equation}\label{1td21}
		P_{2n}(C^\psi_{\bar\beta,2},C)\ge
		\psi(n)\left(\frac1{\sqrt{\pi}}-\sqrt{\frac2{\pi}}\psi(n)\left(\sum\limits_{k=1}^{n-1} \frac{1}{\psi^2(k)}\right)^{\frac12}\right).
	\end{equation}
	The combination of \eqref{1td20} and \eqref{1td21} allows us to write equations
	\begin{equation}\label{1td22}
		P_{2n}(C^\psi_{\bar\beta,2},C)=\psi(n)\left(\frac1{\sqrt{\pi}}+\gamma^{(1)}_n\right),
	\end{equation}
	\begin{equation}\label{1td23}
		P_{2n-1}(C^\psi_{\bar\beta,2},C)=\psi(n)\left(\frac1{\sqrt{\pi}}+\gamma^{(2)}_n\right),
	\end{equation}
	in which for $\gamma^{(i)}_n, i=1,2,$ the following double inequalities hold
	\begin{equation}\label{1td24}
		-\sqrt{\frac2{\pi}}\psi(n)\left(\sum\limits_{k=1}^{n-1} \frac{1}{\psi^2(k)}\right)^{\frac12}\le\gamma^{(i)}_n
		\le
		\frac1{\sqrt{\pi}\psi(n)} \left(\sum\limits_{k=n+1}^\infty \psi^2(k)\right)^{\frac12}.
	\end{equation}
	
	If the condition \eqref{1t3} is satisfied, then by virtue of \eqref{1td22}-\eqref{1td24} the asymptotic equations \eqref{1t4} take place.
	
	Theorem 1 is proved.
\end{proof}

We note that the condition
\begin{equation}\label{25}
	\lim_{n\rightarrow\infty}\frac1{\psi^2(n)}{\sum\limits_{k=n+1}^\infty \psi^2(k)}=0
\end{equation}
is satisfied if $\psi(k)$  satisfies the condition $D_0$:
\begin{equation}\label{26}
	\lim_{k\rightarrow\infty}\frac{\psi(k+1)}{\psi(k)}=0.
\end{equation}

To make sure of this, let us put
\begin{equation*}\label{key}
	\varepsilon_n=\sup_{k\ge n}\frac{\psi(k+1)}{\psi(k)}.
\end{equation*}
By virtue of \eqref{26} $\varepsilon_n\downarrow0$ as $n\rightarrow\infty.$ So, we get
\begin{equation*}\label{key}
	\sum\limits_{k=n+1}^\infty \psi^2(k)=\psi^2(n)\left(\frac{\psi^2(n+1)}{\psi^2(n)}+\frac{\psi^2(n+2)}{\psi^2(n+1)}\frac{\psi^2(n+1)}{\psi^2(n)}+\ldots\right)
	$$
	$$
	\le
	\psi^2(n)\left(\varepsilon^2_n+\varepsilon^4_n+\ldots\right)
	=
	\psi^2(n)\frac{\varepsilon^2_n}{1-\varepsilon^2_n}=o(\psi^2(n)).
\end{equation*}

Let us show that for strictly decreasing sequences $\psi$ the fulfilment of condition $D_0$ of the form \eqref{26}  ensures the truth of the following equality
\begin{equation}\label{27}
	\lim_{n\rightarrow\infty} \psi^2(n) \sum\limits_{k=1}^{n-1} \frac1{\psi^2(k)} =0.
\end{equation}

To do this, we use Stoltz's theorem, according to which the relation \eqref{27} is followed from the following equality
\begin{equation}\label{28}
	\lim_{n\rightarrow\infty} \frac{\sum\limits_{k=1}^{n-1} \frac1{\psi^2(k)}-\sum\limits_{k=1}^{n-2} \frac1{\psi^2(k)}}{\frac1{\psi^2(n)}-\frac1{\psi^2(n-1)}}  =0.
\end{equation}

Since
\begin{equation*}\label{key}
	\frac{\sum\limits_{k=1}^{n-1} \frac1{\psi^2(k)}-\sum\limits_{k=1}^{n-2} \frac1{\psi^2(k)}}{\frac1{\psi^2(n)}-\frac1{\psi^2(n-1)}}
	=
	\frac{\frac1{\psi^2(n-1)}}{\frac1{\psi^2(n)}-\frac1{\psi^2(n-1)}}
	=
	\frac{\frac{\psi^2(n)}{\psi^2(n-1)}}{1-\frac{\psi^2(n)}{\psi^2(n-1)}},
\end{equation*}
then \eqref{28} follows from \eqref{26}.

In view of the above, we have the following statement.

\begin{corollary}
	Let $\bar{\beta}=\{\beta_k\}_{k=1}^\infty, \beta_k\in\mathbb{R},$ and the
	sequence $\psi(k)>0$ is strictly decreasing and satisfies the condition  $D_0$ of the form \eqref{26}. Then as $n\rightarrow\infty$ the asymptotic equalities \eqref{1t4} hold.
\end{corollary}

We give the corollaries of Theorem 1 in some important special cases.

\begin{theorem}\label{theorem2}
	Let $\bar{\beta}=\{\beta_k\}_{k=1}^\infty, \beta_k\in\mathbb{R},$ and $n\in\mathbb{N}.$ Then for all $\displaystyle r\ge\frac{n+1}2$  the following inequalities hold
	$$
	\frac1{\sqrt{\pi}}n^{-r}\left(
	1-\frac{4\left(1-\frac1n\right)^{2r}}{1+4\left(1-\frac1n\right)^{2r}}\right)^{\frac12}
	\le
	P_{2n}(W_{\bar\beta,2}^r,C)
	$$
	\begin{equation}\label{29}		\le
		P_{2n-1}(W_{\bar\beta,2}^r,C)
		\le
		\frac1{\sqrt{\pi}}n^{-r}\left(
		1+\frac{2+\frac1n}{\left(1+\frac1n\right)^{r}}\right)^{\frac12},
	\end{equation}
	where $P_N$ is any of the widths $b_N, d_N, \lambda_N$ or $\pi_N$.
\end{theorem}

\begin{proof}
	Let us put $\psi(k)=k^{-r}, r>1.$ Obviously, the condition \eqref{1t1} is satisfied. Since for $2r\ge n+1, n\in\mathbb{N},$
	$$
	\sum_{k=n+1}^{\infty}\frac1{k^{2r}}<\frac1{(n+1)^{2r}}+\int_{n+1}^{\infty}\frac{dt}{t^{2r}}=\frac1{(n+1)^{2r}}+\frac1{(2r-1)(n+1)^{2r-1}}
	=\frac1{(n+1)^{2r}}\frac{2r+n}{2r-1}
	$$
	\begin{equation}\label{29'}
		\le\frac1{(n+1)^{2r}}\frac{4r-1}{2r-1}
		\le\frac1{(n+1)^{2r}}\left(2+\frac1{2r-1}\right)
		\le
		\frac1{n^{2r}}\frac{2+\frac1{n}}{(1+\frac1n)^{2r}},
	\end{equation}
	then according to the right-hand side of the equality \eqref{1t2} of Theorem 1 we obtain an estimate
	\begin{equation}\label{30}
		P_{2n-1}(W_{\bar\beta,2}^r,C)\le
		\frac1{\sqrt{\pi}}\left(\sum_{k=n}^{\infty}\frac1{k^{2r}}\right)^{\frac12}
		\le 	
		\frac1{\sqrt{\pi}}n^{-r}\left(
		1+\frac{2+\frac1n}{\left(1+\frac1n\right)^{2r}}\right)^{\frac12}.
	\end{equation}
	
	On the other hand, for $\displaystyle r\ge\frac{n+1}2$ and $\psi(k)=k^{-r}$
	$$
	\frac1{\psi^2(n)} +2\sum\limits_{k=1}^{n-1} \frac1{\psi^2(k)}=
	n^{2r}+2\sum\limits_{k=1}^{n-1}k^{2r}
	\le n^{2r}+2\left((n-1)^{2r}+\int\limits_{1}^{n-1}t^{2r}dt\right)
	$$
	$$
	=
	n^{2r}+2\left((n-1)^{2r}+\frac{(n-1)^{2r+1}}{2r+1}-\frac1{2r+1}\right)
	<
	n^{2r}+2\left((n-1)^{2r}+\frac{(n-1)^{2r+1}}{n+2}\right)
	$$
	\begin{equation}\label{31}	
		<
		n^{2r}+4(n-1)^{2r}=n^{2r}\left(1+4\left(1-\frac1n\right)^{2r}\right).
	\end{equation}
	
	By virtue of the left part of the inequality \eqref{1t2} of Theorem 1 and the formula \eqref{31} we get an estimate
	\begin{equation}\label{32}
		P_{2n}(W_{\bar\beta,2}^r,C)
		\ge	
		\frac1{\sqrt{\pi}}n^{-r}\left(
		\frac{1}{1+4\left(1-\frac1n\right)^{2r}}\right)^{\frac12}
		=	\frac1{\sqrt{\pi}}n^{-r}\left(
		1-\frac{4\left(1-\frac1n\right)^{2r}}{1+4\left(1-\frac1n\right)^{2r}}\right)^{\frac12}.
	\end{equation}
	Combining the estimates \eqref{30} and \eqref{32} we obtain \eqref{29}.
	Theorem 2 is proved.
\end{proof}

Note that if a condition
\begin{equation}\label{33}
	\lim\limits_{n\rightarrow\infty}\frac rn=\infty
\end{equation}
is satisfied, then for $\psi(k)=k^{-r}$ the condition \eqref{26} is also satisfied because
\begin{equation*}\label{key}
	\frac{\psi(k+1)}{\psi(k)}=\left(\frac k{k+1}\right)^r=\left(1+\frac1k\right)^{-r}=\left(\left(1+\frac1k\right)^{k+1}\right)^{-\frac r{k+1}}\le e^{-\frac r{k+1}}\rightarrow0,\ k\rightarrow\infty.
\end{equation*}

Taking the limit, as $n\rightarrow\infty$, in the relations \eqref{29}, we obtain the following statement.

\begin{theorem}\label{theorem3}
	Let $\bar{\beta}=\{\beta_k\}_{k=1}^\infty, \beta_k\in\mathbb{R},$ $n\in\mathbb{N},$ and the condition \eqref{33} is satisfied.
	Then the following asymptotic equalities hold
	\begin{equation}\label{34}
		\left.\begin{array}{l}
			P_{2n}(W_{\bar\beta,2}^r,C)\\
			P_{2n-1}(W_{\bar\beta,2}^r,C)
		\end{array}\right\}=
		n^{-r}\left(\frac1{\sqrt{\pi}}+
		{\mathcal O}(1)\left(1+\frac1{n}\right)^{-r}\right),
	\end{equation}
	where $P_N$ is any of the widths $b_N, d_N, \lambda_N$ or $\pi_N$,  and ${\mathcal O}(1)$ are the quantities uniformly bounded in all parameters.
\end{theorem}

Note also that the equalities \eqref{34} are easy obtained from the formula \eqref{1t4} and estimates \eqref{29'} and \eqref{31}.

\begin{theorem}\label{theorem4}
	Let $\bar{\beta}=\{\beta_k\}_{k=1}^\infty, \beta_k\in\mathbb{R}$, $\alpha>0,$ $r>1,$ $n\in\mathbb{N}$ and be such that
	\begin{equation}\label{4t1}
		(n-1)^r>\frac1\alpha,
	\end{equation}
	then the following inequalities hold
	$$
	\frac1{\sqrt{\pi}}e^{{-}\alpha n^r}\left(
	1{-}\frac{2\gamma_{\alpha,r,n}e^{-2\alpha r (n-1)^{r-1}}}{1+2\gamma_{\alpha,r,n}e^{-2\alpha r (n-1)^{r-1}}}
	\right)^{\frac12}
	\le
	P_{2n}(C^{\alpha,r}_{\bar\beta,2},C)
	$$
	\begin{equation}\label{'29}
		\le
		P_{2n-1}(C^{\alpha,r}_{\bar\beta,2},C)
		\le	
		\frac1{\sqrt{\pi}}e^{-\alpha n^r}\left(1+e^{-2\alpha r n^{r-1}}\left(1+\frac1{2\alpha r n^{r-1}}\right)\right)^{\frac12},
	\end{equation}
	where $P_N$ is any of the widths $b_N, d_N, \lambda_N$ or $\pi_N$ and 
	\begin{equation}\label{4t4}
		\gamma_{\alpha,r,n}=\left(1+\frac{1}{\alpha r (n-1)^{r-1}}+e^{-2\alpha (n-1)^{r}}\max \left\{{e^{4\alpha}},\frac{e^2}{\alpha^{1{+}1/r}}\right\}\right).
	\end{equation}
\end{theorem}

\begin{proof}
	First of all, note that if $\alpha>0,$ $r>1,$ $n\in\mathbb{N}$ and satisfy the condition \eqref{4t1}, then for a quantity of the form
	\begin{equation}\label{p4t1}
		I_{n-1}:=\int\limits_{1}^{n-1}e^{2\alpha t^{r}}dt,\quad \alpha>0,\ r>1,
	\end{equation}
	the following inequality holds
	\begin{equation}\label{p4t2}
		I_{n-1}\le\frac{e^{2\alpha (n-1)^{r}}}{\alpha r (n-1)^{r-1}}+\max \left\{{e^{4\alpha}},\frac{e^2}{\alpha^{1{+}1/r}}\right\}.
	\end{equation}
	Indeed, integrating by parts we have
	\begin{equation}\label{p4t3}
		I_{n-1}=\frac1{2\alpha r}\int\limits_{1}^{n-1}  t^{1-r}de^{2\alpha t^{r}}=\frac1{2\alpha r}\left(\frac{e^{2\alpha (n-1)^{r}}}{(n-1)^{r-1}}-e^{2\alpha}\right)+\frac{r-1}{2\alpha r}\int\limits_{1}^{n-1}\frac{e^{2\alpha t^{r}}}{t^r}dt.
	\end{equation}
	For $0<2\alpha<1$, taking into account \eqref{4t1}, we obtain
	\begin{equation}\label{p4t4}
		\int\limits_{1}^{n-1}\frac{e^{2\alpha t^{r}}}{t^r}dt
		=
		\int\limits_{1}^{\alpha^{-1/r}}\frac{e^{2\alpha t^{r}}}{t^r}dt+\int\limits_{\alpha^{-1/r}}^{n-1}\frac{e^{2\alpha t^{r}}}{t^r}dt
		<\frac{e^2}{\alpha^{1/r}} +\alpha\int\limits_{\alpha^{-1/r}}^{n-1}e^{2\alpha t^{r}}dt\le \frac{e^2}{\alpha^{1/r}}+\alpha I_{n-1}.
	\end{equation}
	From \eqref{p4t3} and \eqref{p4t4} under condition  \eqref{4t1} 
	we get
	$$
	I_{n-1}<\frac{e^{2\alpha (n-1)^{r}}}{2\alpha r(n-1)^{r-1}}+\frac {e^2}{2\alpha^{1+1/r}}+\frac12I_{n-1}
	$$
	and
	\begin{equation}\label{p4t5}
		I_{n-1}<\frac{e^{2\alpha (n-1)^{r}}}{\alpha r(n-1)^{r-1}}+\frac {e^2}{\alpha^{1+1/r}},\quad 0<2\alpha<1.
	\end{equation}
	
	For $2\alpha\ge1$ we obtain
	$$
	\int\limits_1^{n-1}\frac{e^{2\alpha t^{r}}}{t^r}dt= \int\limits_{1}^{2^{1/r}}\frac{e^{2\alpha t^{r}}}{t^r}dt+\int\limits_{2^{1/r}}^{n-1}\frac{e^{2\alpha t^{r}}}{t^r}dt\le\int\limits_{1}^{2^{1/r}}{e^{2\alpha t^{r}}}dt+\frac12\int\limits_{2^{1/r}}^{n-1}{e^{2\alpha t^{r}}}dt=
	$$
	\begin{equation}\label{p4t6}
		=\frac12\int\limits_{1}^{2^{1/r}}{e^{2\alpha t^{r}}}dt+\frac12I_{n-1}<
		\frac{2^{1/r}-1}{2}e^{4\alpha}+\frac12I_{n-1}<\frac{e^{4\alpha}}2+\frac12I_{n-1}.
	\end{equation}
	From \eqref{p4t3} and \eqref{p4t6} under condition  \eqref{4t1} we have
	$$
	I_{n-1}<\frac{e^{2\alpha (n-1)^{r}}}{2\alpha r(n-1)^{r-1}}+\frac {e^{4\alpha}}{4\alpha}+\frac1{4\alpha}I_{n-1}<\frac{e^{2\alpha (n-1)^{r}}}{2\alpha r(n-1)^{r-1}}+\frac {e^{4\alpha}}{2}+\frac12I_{n-1}
	$$
	and
	\begin{equation}\label{p4t7}
		I_{n-1}<\frac{e^{2\alpha (n-1)^{r}}}{\alpha r(n-1)^{r-1}}+ {e^{4\alpha}},\quad 2\alpha\ge1.
	\end{equation}
	
	The inequality \eqref{p4t2} follows from \eqref{p4t5} and \eqref{p4t7}.
	
	For  $\psi(k)=e^{-\alpha k^{r}},\ \alpha>0, \ r>1, $ under condition  \eqref{4t1}, taking into account that for $r>1$ and $ n\in\mathbb{N}$
	\begin{equation*}\label{p4t8}
		\left(1+\frac 1n\right)^r-1>\frac {r}n,
	\end{equation*}
	we obtain 
	$$
	\psi^2(n)\sum\limits_{k=1}^{n-1}\frac1{\psi^2(k)}=
	e^{-2\alpha n^{r}}\sum\limits_{k=1}^{n-1}	e^{2\alpha k^{r}}
	\le e^{-\alpha n^{r}}\left( e^{2\alpha (n-1)^{r}}+\int\limits_{1}^{n-1}e^{2\alpha t^{r}}dt\right)
	$$$$
	\le e^{-2\alpha n^{r}} \left(e^{2\alpha (n-1)^{r}}\left(1+\frac{1}{\alpha r (n-1)^{r-1}}+e^{-2\alpha (n-1)^{r}}\max \left\{{e^{4\alpha}},\frac{e^2}{\alpha^{1{+}1/r}}\right\}\right)\right)
	$$
	\begin{equation}\label{p4t9}
		\le
		e^{-2\alpha (n-1)^{r}((1+\frac1{n-1})^r-1)} \gamma_{\alpha,r,n}
		\le
		\gamma_{\alpha,r,n}e^{-2\alpha r (n-1)^{r-1}}.
	\end{equation}
	
	Thus, by virtue of the left part of the inequality \eqref{1t2} of Theorem 1 and  \eqref{p4t9},  we obtain a required estimate from below for widths $P_{2n}(C^{\alpha,r}_{\bar\beta,2},C), \ \alpha>0, \ r>1,$   under condition  \eqref{4t1}
	$$
	P_{2n}(C^{\alpha,r}_{\bar\beta,2},C)\ge 
	\frac1{\sqrt{\pi}}\left(e^{2\alpha n^{r}}+2\sum\limits_{k=1}^{n-1}	e^{2\alpha k^{r}}\right)^{-\frac12}
	$$
	\begin{equation}\label{p4t10}
		\ge\frac1{\sqrt{\pi}}e^{{-}\alpha n^r}\left(
		1{-}\frac{2\gamma_{\alpha,r,n}e^{-2\alpha r (n-1)^{r-1}}}{1+2\gamma_{\alpha,r,n}e^{-2\alpha r (n-1)^{r-1}}}
		\right)^{\frac12}.
	\end{equation}
	
	As was shown in \cite[P. 163-164]{Stepanets1987}
	\begin{equation}\label{24}
		\sum_{k=n+1}^{\infty}e^{-\alpha k^r}<e^{-\alpha n^r}\left(1+\frac1{\alpha r n^{r-1}}\right)e^{-\alpha r n^{r-1}}, \quad r>1, \alpha>0, n\in\mathbb{N}.
	\end{equation}
	Therefore,
	\begin{equation}\label{p4t11}
		\frac1{\psi^2(n)}\sum\limits_{k=n+1}^\infty\psi^2(k)=e^{2\alpha n^r}\sum\limits_{k=n+1}^\infty e^{-2\alpha k^r}
		%\le	
		%e^{\alpha n^r}\sum_{k=n+1}^{\infty}e^{-\alpha k^r}
		<\left(1+\frac1{2\alpha r n^{r-1}}\right)e^{-2\alpha r n^{r-1}}.
	\end{equation}
	
	Thus,  by virtue of the right part of the inequality \eqref{1t2} of Theorem 1 and the formula \eqref{p4t11} we get an estimate
	\begin{equation}\label{p4t12}
		P_{2n}(W_{\bar\beta,2}^r,C)
		\le	
		\frac1{\sqrt{\pi}}e^{-\alpha n^r}\left(1+e^{-2\alpha r n^{r-1}}\left(1+\frac1{2\alpha r n^{r-1}}\right)\right)^{\frac12}.
	\end{equation}
	Combining the estimates \eqref{p4t10} and \eqref{p4t12} we obtain \eqref{'29}.
	Theorem 4 is proved.
\end{proof}

Taking the limit, as $n\rightarrow\infty$, in the relations \eqref{'29}, we obtain the following statement.

\begin{theorem}\label{theorem5}
	Let $\bar{\beta}=\{\beta_k\}_{k=1}^\infty, \beta_k\in\mathbb{R}$, $\alpha>0,$ $r>1,$ $n\in\mathbb{N}$ and the condition \eqref{4t1} is satisfied.  Then as $n\rightarrow\infty$ the following asymptotic equalities hold	
	\begin{equation}\label{5t1}
		\left.\begin{array}{l}
			P_{2n}(C^{\alpha,r}_{\bar\beta,2},C)\\
			P_{2n-1}(C^{\alpha,r}_{\bar\beta,2},C)
		\end{array}\right\}=
		e^{-\alpha n^r}\left(\frac1{\sqrt{\pi}}
		+ {\mathcal O}(1)\gamma_{\alpha,r,n}	e^{-\alpha r (n-1)^{r-1}}\right),
	\end{equation}
	where $P_N$ is any of the widths $b_N, d_N, \lambda_N$ or $\pi_N$ and $\gamma_{\alpha,r,n}$  is defined by \eqref{4t4} and ${\mathcal O}(1)$ are the quantities uniformly bounded in all parameters.
\end{theorem}

Note that the Theorem \ref{theorem5} complements the results of the works \cite{Serdyuk1999,Serdyuk_Bodenchuk2013JAT,Serdyuk_Sokolenko2011JAT,Shevaldin_1992,Stepanets_Serdyuk1995UMZh8},   which contain exact estimates for the widths of the classes of convolutions  with classical or generalized Poisson kernels. As it follows from the proofs of Theorems 4 and 5 the asymptotic equalities for widths in \eqref{5t1} are realized by trigonometric Fourier sums. The asymptotic equalities for deviations of Fourier sums on classes of generalized Poisson integrals $C^{\alpha,r}_{\beta,p}$ in the uniform metric are seen, for example, in \cite{Serdyuk_2005_8,Serdyuk_Stepanyuk2019,Stepanets2005} and others.

\bigskip

CONTACT INFORMATION

\medskip
A.S.~Serdyuk\\01024, Ukraine, Kyiv-4, 3, Tereschenkivska st.\\
serdyuk@imath.kiev.ua

\medskip
I.V.~Sokolenko\\01024, Ukraine, Kyiv-4, 3, Tereschenkivska st.\\
sokol@imath.kiev.ua
\end{document}